\documentclass[11pt]{amsart}
\usepackage{amssymb,amsfonts,amsthm,amsmath}
\usepackage{mathrsfs}
\usepackage[english]{babel}  

\newtheorem{lemma}{Lemma}[section]
\newtheorem{corollary}[lemma]{Corollary}
\newtheorem{definition}[lemma]{Definition}
\newtheorem{proposition}[lemma]{Proposition}
\newtheorem{remark}[lemma]{Remark}
\newtheorem{theorem}[lemma]{Theorem}

\newcommand{\field}[1]{\mathbb{#1}}
\newcommand{\B}{{\bf B}}

\newcommand{\C}{\field{C}}
\newcommand{\FH}{\field{H}}
\newcommand{\FP}{\field{P}}
\newcommand{\G}{{\bf G}}
\newcommand{\GR}{{\bf G}_\R} 
\newcommand{\K}{\field{K}}
\newcommand{\N}{\field{N}}
\newcommand{\R}{\field{R}}
\newcommand{\Sr}{{\bf S}}
\newcommand{\Srn}{{\bf S}^{n-1}}
\newcommand{\Om}{\Omega}
\newcommand{\Rn}{\R^n}

\newcommand{\Cal}[1]{\mathcal{#1}}
\newcommand{\Scr}[1]{\mathscr{#1}}

\newcommand{\nf}{\nabla f}

\newcommand{\bh}{{\bf h}}
\newcommand{\bu}{{\bf u}}

\newcommand{\taug}{{\tau}_g}
\newcommand{\tauf}{{\tau}_f}
\newcommand{\tf}{{t}_f}
\newcommand{\tg}{{t}_g}
\newcommand{\dr}{\partial_r}
\newcommand{\ft}{\tilde{f}}
\newcommand{\np}{\nabla^{'}}

\newcommand{\ud}{{\rm{d}}}

\newcommand{\vol}{{\rm{vol}}}
\newcommand{\clos}{{\rm{\bf clos}}}
\newcommand{\crit}{{\rm{\bf crit}}}
\newcommand{\fl}{\rightarrow}

\newcommand{\bs}{ {\tiny $\blacksquare$} \\}
\numberwithin{equation}{section}

\newenvironment{myproof}{\noindent{\it Proof.}
\setlength{\parindent}{0mm}}
{$\hfill$ \bs}

\begin{document}
\title[Tame Functions with strongly isolated singularities at infinity]{Tame Functions with strongly
isolated singularities at infinity: a tame version of a Parusi\'nski's Theorem}

\author[]{Vincent Grandjean}

\address{{\it Permanent Address:} V. Grandjean, Department of Computer Science,
University of Bath, BATH BA2 7AY, England,(United Kingdom)}

\address{{\it Current Address:} V. Grandjean, Fakul\"at V, Institut f\"ur Mathematik Carl
von Ossietzky Universit\"at, Oldenburg, 26111 Oldenburg i.O. (Germany)}

\email{cssvg@bath.ac.uk}



\begin{abstract}
Let $f:\Rn \mapsto \R$ be a definable function, enough
differentiable. Under the condition of having strongly isolated singularities at infinity at a regular
value $c$, we give a sufficient condition expressed in terms of the total absolute curvature function to
ensure the local triviality of $f$ over a neighbourhood of $c$ and doing so providing the tame version of
Parusi\'nski's Theorem on complex polynomials with isolated singularities at infinity.
\end{abstract}
\maketitle

\section{Introduction}
The fundamental result of Thom about the finiteness of the topological types of a given polynomial
function \cite{Th}, has led to some understanding of the geometry of the foliation by the level
of a given tame function $f:U \subset \K^n \mapsto \K$ nearby a generalised critical value, value likely to
be a bifurcation value, that is at which the topology of the fibres is not locally constant. Rather early
one noticed that a bifurcation value could be a regular value, as already suggested by the properness
condition in Erehsmann's Theorem to ensure the local triviality of a submersion. For a decade or so, there
were no "effective" criterion to describe these regular bifurcation values, or at least a finite subset
of $\K$ that would contain them. Then came some sufficient conditions to trivialise the given function
over a regular value (see \cite{Ph}, \cite{Br}, \cite{HL}). This led rapidly to the notion of asymptotic
critical value (or generalised critical value), requiring, similarly to the vanishing of the gradient
at a critical point, that the gradient vector field is asymptotically small along a sequence going to
the boundary of the domain. To be more precise there exists a sequence $x$ going to "infinity" along
which $f(x) \fl c \in \K$ and $|x|\cdot |\nf (x)| \fl 0$.

\medskip
For a tame function defined on $\K^n$ the boundary must be understood as infinity, more precisely the hyperplane 
at infinity of the usual projective compactification of $\K^n$. In this context, a
regular bifurcation value seems to be like a critical value coming from some sort of singular phenomenon
(to be fully understood) lying on the boundary. Then it was proved (\cite{Ph}, \cite{HL}, \cite{Pa1},
\cite{Ti}, \cite{LZ}, \cite{dac}, etc...) that any regular bifurcation value must be an asymptotic critical
value and there are finitely many such asymptotic critical values.

H\`a and L\^e proved in \cite{HL} that the triviality of the complex plane polynomial function $f$ over
a neighbourhood of a given value $c$ was equivalent to the constancy of the Euler Characteristic of the
fibres in a neighbourhood of the value $c$. This result was later generalised by Parusi\'nski to the case
of complex polynomials with isolated singularities at infinity \cite{Pa1}, which he also proved to be
equivalent to requiring that the regular value $c$ is not an asymptotic critical value. Parusi\'nski
\cite{Pa2} also explained that being an asymptotic critical value of $f$ is equivalent to the failure
of a certain stratifying condition on the projective closure of the graph of $f$, expressed as a
property of the relative conormal space of $f$, what Tib\u{a}r has called $t$-isolated singularities
(\cite{Ti}). Despite the same property holds true for real polynomials (\cite{Ti}), the phenomena
occurring in the real domain are of a much less rigid kind. There are already counter-example of
H\`a-L\^e's result as noticed by Tib\u{a}r and Zaharia \cite{TZ}. They nevertheless provided necessary
and sufficient conditions for a real plane polynomial function to be locally trivial over a neighbourhood
of a value $c$. Later, in \cite{CP}, Coste and de la Puente proved an equivalent version of
Tib\u{a}r-Zaharia's result in terms of polar curves, that is involving the relative conormal geometry
of the function at infinity nearby the given level $c$. In the world of real tame functions (to be
understood as "globally" definable in an o-minimal structure expanding the ordered field of real
numbers), the hope to find necessary and sufficient conditions for a regular value to be a bifurcation value
is much harder to hold ! Different sufficient conditions were provided to guarantee the local triviality
at infinity over a neighbourhood of a regular value (see \cite{LZ}, \cite{Ti}, \cite{TZ}, \cite{dac},
\cite{DDVG1}, \cite{DDVG2}, \cite{DDVG3}). There are unfortunately not all comparable, but they all exhale
a similar flavour: The naive belief that too much bending (curvature and so a possible lack of transversality
to "spheres") is an obstruction to trivialisation. Thus the understanding of the relative conormal geometry
at infinity nearby a regular level $c$ we are interested in is an important aspect to explore in order to
deciding whether the value $c$ is a bifurcation value or not.

\medskip
The aim of this paper is to provide a real version in the (globally) definable setting of Parusi\'nski's result.

We propose a condition on the fibre, at a regular value $c$, of the relative conormal space of a tame
function $f$, that we call {\bf SISI} at $c$ (shortening for strongly isolated singularities at infinity, see
Definition \ref{definTITFSISI1}). This condition geometrically means:
For a given value $c$, there are at most finitely many points at infinity such that any limit 
of tangent hyperplanes, along any sequence going to infinity with limit of secants the given point at infinity 
and along which the function tends to $c$, may not be orthogonal to the line direction corresponding to this point.

This property, about the asymptotic behaviour of limits of tangent hyperplanes to
the fibres of $f$ when getting closer and closer to the level $c$ in a neighbourhood of infinity, when
combined with a property on the total absolute curvature of the function $f$, ensures the triviality of $f$
nearby $c$ (Theorem \ref{thmTITFSISI1}).
This condition is satisfied for a real polynomial having isolated
singularities at infinity as defined by Parusi\'nski. We have stated our result in terms of the continuity
of the total absolute curvature function of $f$ which is the same sort of condition of having the
generic polar curves empty at infinity in a neighbourhood of the level $c$.

\smallskip The paper is organised as follows:

\smallskip We begin with Section \ref{sectionNC} in which we explain some notations and some conventions
that will be later used in the paper.

\smallskip Section \ref{sectionGMTF}, \ref{sectionRCGITF} and \ref{sectionACVBV} are reminders of the
definitions and of some of the elementary properties of the key objects we work with, such as the total absolute
curvature function, the relative conormal space and the notion of asymptotic critical value, we want to
focus on in Sections \ref{sectionTITFSISI} and \ref{sectionRPCP}.

\smallskip Condition {\bf SISI} is defined in Section \ref{sectionTITFSISI}, where we state and proof our
main result:

\smallskip 
\noindent 
{\bf Theorem \ref{thmTITFSISI1}}. {\em Let $f:\Rn\mapsto\R$ be a $C^l$ definable function with
$l\geqslant 2$. Assume that
the function $f$ satisfies condition {\bf SISI} at a regular value $c$. \\
If the function  $t \mapsto |K|(t)$ is continuous at $c$, we trivialise $f$ over a neighbourhood of $c$
by means of the flow of a $C^{l-1}$ vector field. So $c \notin B(f)$.}

\smallskip In Section \ref{sectionRPCP}, we compare condition {\bf SISI} for a real polynomial function and
Parusi\'nski's notion of {\it isolated singularities at infinity} (defined for complex polynomials). Our
other main result is

\smallskip 
\noindent
{\bf Proposition \ref{propRPCP2}}. {\em If $f:\Rn\mapsto\R$ is a real polynomial with isolated
singularities at infinity, then condition {\bf SISI} at any regular value $c$ is satisfied}.

\smallskip We finish with some final remarks and comments in Section \ref{sectionCR} to explain that our
result is really Parusi\'nski real counter-part and, unfortunately, nothing much better is to be expected for
real polynomials that what we do not already have with this level of generality.
%
%
%
%
%
%
%
%
%
%
%
%
%
\section{Notation - convention}\label{sectionNC}
Let $\Rn$ be the real $n$-dimensional affine space endowed with its Euclidean metric.  The scalar product
will be denoted by $\langle \cdot , \cdot \rangle$.

Let $\FP_\K^n$ be the projectivised space of the $\K$-vector space $\K^n$. This notation will be
exclusively used to mean "the projective" compactification of $\K^n$, where $\K$ either stands for $\R$ or
$\C$.

Let $\B_R^n$ be the open ball of $\R^n$ centred at the origin and of radius $R >0$.

Let $\Sr_R^{n-1}$ be the $(n-1)$-sphere centred at the origin and of radius $R>0$.

Let $\Sr^{n-1}$ be unit ball of $\R^n$.

Let $\bh$ be a $\K$-vector subspace of dimension $q$ of a $\K$-vector space $E$.
Let $\G_\K (p,\bh)$ be the Grassmann manifold of the $p$-dimensional $\K$-vector subspaces of $\bh$. When $\bh =
E = \K^q$, we will only write $\G_\K (p,q)$.

\medskip
Let us recall briefly what an o-minimal structure is.

\medskip
An {\it o-minimal structure $\Cal{M}$ expanding the ordered field of real numbers} is a collection
$(\Cal{M}_p)_{p \in \N}$, where $\Cal{M}_p$ is a set of subsets of $\R^p$ satisfying the following axioms

\smallskip \noindent
1) For each $p\in \N$, $\Cal{M}_p$ is a boolean subalgebra of subsets of $\R^p$. \\
2) If $A \in \Cal{M}_p$ and $B \in \Cal{M}_q$, then $A \times B \in \Cal{M}_{p+q}$. \\
3) If $\pi: \R^{p+1} \mapsto \R^p$, is the projection on the first $p$ factors, given any $A \in
\Cal{M}_{p+1}$, $\pi (A) \in \Cal{M}_p$.  \\
4) The algebraic subsets of $\R^p$ belongs to $\Cal{M}_p$. \\
5) $\Cal{M}_1$ consists exactly of the finite unions of points and intervals.

So the smallest o-minimal structure is the structure of the semi-algebraic subsets.

\medskip
Assume that such an o-minimal structure $\Cal{M}$ is given for the rest of this article.

\medskip
A subset $A$ of $\R^p$ is a {\it definable subset} (in the given o-minimal structure) of $\R^p$, if $A \in
\Cal{M}_p$.

A subset $B$ of $\FP_\R^p$ is said to be {\it (globally) definable} if its trace in each affine chart is 
a definable
subset of this chart.

For mappings we  slightly restrict the usual definition of a definable mapping. A mapping $g: X \mapsto Y$,
where $X \subset \FP_R^p$ and $Y \subset \R^q$, is a {\it definable mapping} (or just definable, for short)
if the intersection of the closure of its graph in $\FP_\R^p \times\R^q$ with $\FP_\R^p \times \B$ is a
definable subset of $\FP_\R^p \times\R^q$ for any ball $\B \subset \R^q$.

The reader may refer to \cite{Cos,vdD1,vdDM} to learn more about the properties of definable subsets and
definable mappings.
%
%
%

\medskip
Let $S$ be a $C^1$ definable submanifold of $\FP_\R^n$. Let $g : S \mapsto \R^q$ be a $C^1$ definable
mapping. The critical set of $g$ is denoted by $\crit (g)$.

By abuse of language, we will talk about the rank of the mapping $g$
at a point $x_0$ to mean the rank of the differential $\ud_{x_0}g$. We will also talk about the rank of $g$
to mean the maximal rank of the differentials $\ud_xg$, $x\in S$.

\begin{remark}
In this paper, we will {\rm always} use the adjective {\rm definable} for a subset of an affine space to
mean {\rm definable} in the projective compactification of the ambient affine space.
\end{remark}

Let $\varphi$ and $\psi$ be two germs at the origin (resp. at infinity) of single real variable functions. We
write $\varphi \sim \psi$ to mean that the ratio $\varphi/\psi$ has a non zero finite limit at the origin 
(resp. at infinity). We write $\varphi \simeq \psi$ when the limit of $\varphi/\psi$ at the origin (resp. 
at infinity) is $1$.
We will write $\psi = o(\varphi)$ to mean $\psi/\varphi \fl 0$ at the origin (resp. at infinity).
%
%
%
%
%
%
%
%
%
%
%
%
%
%
%
\section{Gauss Map of a tame function}\label{sectionGMTF}

Let $f:\Rn \mapsto \R$ be a $C^l$ definable function, with $l \geqslant 2$.

Let $\crit (f)$ be the critical set of the function $f$ and let $K_0 (f)$ be the set of its critical 
values, that is $K_0 (f) = f(\crit (f))$, that we recall is finite.

For each $t$, let $F_t$ be the level $f^{-1} (t)$.

\medskip
The Gauss map of the function $f$ is the mapping defined as follows:

\medskip
\begin{tabular}{rcccc}
\hspace{3cm} $\nu_f$ &  :& $\R^n \setminus \crit (f)$ &  $\mapsto$&
$\Srn$ \\
& & $x$ & $\mapsto$ & $\displaystyle{
\frac{\nabla f (x)}{|\nabla f (x)|}}$
\end{tabular}

\medskip
It is a definable mapping that is $C^{l-1}$.
Thus the set of its critical values $\nu_f(\crit (\nu_f))$ is
a definable subset of dimension at most $n-2$.

\medskip
For each regular value $t$, let $\nu_t$ be the restriction of $\nu_f$ to $F_t$. So it is also a Gauss 
map on each connected components of $F_t$ providing each component with an orientation that is compatible
with the transverse structure of the foliation of $\Rn \setminus \crit (f)$ by the levels of the
function $f$. Note also that $\crit (\nu_f) \cap F_t = \crit (\nu_t)$.

\medskip
For a given $x \in F_t$, let $k_t (x)$ be the Gaussian curvature of $F_t$ at $x$, namely
$k_t (x) = {\rm det} (\ud_x \nu_t)$.

\medskip
Let $\ud v_{n-1}$ be the $(n-1)$-dimensional Hausdorff measure of $\R^n$.

\begin{definition} Let $t$ be a regular value of the function $f$.

(i) The total absolute curvature of the level $F_t$ is
\begin{center}
\vspace{4pt}
$|K|(t) = \int_{F_t} |k_t (x)| \ud v_{n-1} (x)$
\vspace{4pt}
\end{center}

(ii) The total  curvature of the level $F_t$ is
\begin{center}
\vspace{4pt}
$K(t) = \int_{F_t} k_t (x) \ud v_{n-1} (x)$
\vspace{4pt}
\end{center}
\end{definition}

\medskip
Let us say few words about these total curvatures.
First they are well defined as it will appeared below.

Let us define $\Psi_f : \Rn \setminus (\crit (f) \cup\crit (\nu_f)) \mapsto \Srn \times \R$ such that  $x
\mapsto \Psi_f (x) :=(\nu_f(x),f(x))$.

Let $\widetilde{\Cal{U}}$ be the image of $\Psi_f$. It is an open definable subset since $\Psi_f$ is a
local diffeomorphism at each of its point. The subset $\widetilde{\Cal{U}}$  is a finite disjoint union
$\sqcup_m \widetilde{\Cal{U}}_m$, where $\widetilde{\Cal{U}}_m = \{(u,t) : \# \Psi_f^{-1} (u,t) = m\}$. For
any regular value $t$, let $\Cal{U}_t := \{u \in \Srn: (u,t) : \in \widetilde{\Cal{U}}\}$. 
It is an open definable subset of
$\Srn$. Let  $(\Cal{U}_{i,t})_{i=1,\ldots,q_t}$ be the set of connected components of $\Cal{U}_t$. For each
$i = 1,\ldots,q_t$, let $m_{i,t}$ be $\# \Psi_f^{-1}(u,t)$, for any $u \in \Cal{U}_{i,t}$.

\smallskip
From Gabrielov's uniformity principle there exists a positive integer $N_f$ such that for each $(u,t)
\in \widetilde{\Cal{U}}$, we deduce $\# (\nu_t^{-1} (u) \cap (F_t \setminus \crit (\nu_f))) \leqslant N_f$. 
So we find
\begin{center}
\vspace{6pt} $|K|(t) =\sum_{i=1}^{q_t} m_{i,t}\vol_{n-1} (\Cal{U}_{i,t})$, and $K(t) =
\sum_{i=1}^{q_t}\delta_{i,t} \vol_{n-1} (\Cal{U}_{m,i}^t)$, 
\vspace{6pt}
\end{center}
where $\delta_{i,t}$ is the degree of the mapping $\nu_t$ at any $x \in \nu_t^{-1} (u)$,
$u \in \Cal{U}_{i,t}$.

Thus we have defined two functions $|K| : \R \setminus K_0 (f) \mapsto \R$, $t \mapsto |K| (t)$, the total
absolute curvature function and $K : \R \setminus K_0 (f) \mapsto \R$, $t \mapsto K (t)$, the total
curvature function.

It is a matter of interest to know more about the regularity properties of these two functions, since
they are closely linked to the topology of the levels $F_t$, as the usual Gauss-Bonnet-Chern Theorem suggests
in the compact connected odd-dimensional case.

In this general setting little is known about such functions, nevertheless we know that

\begin{theorem} [\cite{Gr2}] Let $f: \Rn \mapsto \R$ be a $C^l$ definable function, with $l \geqslant 2$.

(i) The function $t \mapsto |K| (t)$ has at most finitely many discontinuities.

(ii) If the function $t \mapsto |K| (t)$ is continuous at a regular value $c$, so is $t \mapsto K (t)$.
\end{theorem}

Obviously if the Gauss map $\nu_f$ is degenerate, that is of rank at most $n-2$, the former
results are without interest since both total curvature functions are the null function.

Since there are finitely many values at which $|K|$ may not
be continuous and finitely many values at which the topology of the fibres of $f$ is not locally constant,
is there a link between these two set of values ? \\
We will see in the next sections that with additional hypotheses there are such relations.

To finish this section let us state the following result that will be important in Section \ref{sectionTITFSISI}.

\begin{proposition}[{\cite[Corollary 6.3]{Gr2}}]\label{propGMTF1}
Let $c$ be a regular value at which the function $|K|$ is not continuous. There exists an open subset
$U \subset \Srn$, such that for any $u \in U$, there exists a connected component $\Gamma$ of
$\Psi_f^{-1}(\{u\}\times \R)$, such that $\Gamma \cap f^{-1}(c)$ is not empty and one of the
two situations below happens:

(i) If $c$ is the infimum of $f$ along $\Gamma$, for any
$\varepsilon >0$ small enough, $\Gamma \cap
f^{-1}(]c,c+\varepsilon[)$ is not bounded.

(ii) If $c$ is the supremum of $f$ along $\Gamma$, for any
$\varepsilon >0$ small enough, $\Gamma \cap
f^{-1}(]c-\varepsilon,c[)$ is not bounded.
\end{proposition}

Note that obviously the oriented polar curve $\Psi_f^{-1}(\{u\}\times \R)$ is $C^{l-1}$ and has at most $N_f$
connected components lying in $\Rn \setminus (\crit(f) \cup \crit (\nu_f))$.
%
%
%
%
%
%
%
%
%
%
%
%
\section{Relative conormal geometry at infinity of a tame function}\label{sectionRCGITF}

Let $\FH_\K^\infty := \FP_\K^n \setminus \K^n$ be the {\it hyperplane at infinity}.

\begin{definition}
Let $g : S \mapsto \R$ be a $C^l$ definable mapping from a submanifold $S\subset\Rn$. The relative
projective conormal bundle of the function $g$ is the subset $\Scr{X}_g$ of $\FP_\R^n \times \GR
(n-1,n)\times \R$ defined as the closure of
\begin{center}
\vspace{4pt}
$\{(x,\bh,t) \in (S \setminus \crit (g)) \times \GR (n-1,n) \times \R : T_x g \subset \bh, t = g(x)\}$,
\vspace{4pt}
\end{center}
where $T_x g = T_x (g^{-1} (g(x)))$.
\end{definition}

The subset $\Scr{X}_g$ is a closed definable subset of $\FP_\R^n \times \GR (n-1,n)\times \R$ of dimension
$n$. Note that $\Scr{X}_g \cap ((S\setminus \crit (g)) \times \GR (n-1,n) \times \R)$ is a $C^{l-1}$
submanifold of $\FP_\R^n \times \GR (n-1,n) \times \R$ since the $x \mapsto T_x g \in \G_\R (\dim S,n)$
is just the projective Gauss map.

\medskip
Let $(\pi_g,\taug,\tg) :\Scr{X}_g \mapsto\FP_\R^n \times \GR (n-1,n) \times \R$ be the restriction of
projections on the respective factors of $\FP_\R^n \times \GR (n-1,n) \times \R $, that is, 
\begin{center}
\vspace{4pt}
$\pi_g (x,t,\bh)
= x$, $\taug (x,t,\bh) =\bh$ and $\tg (x,t,\bh) = t$. 
\vspace{4pt}
\end{center}
Those maps are definable, and $C^{l-1}$ on $\Scr{X}_g \cap ((S \setminus
\crit (g)) \times \GR (n-1,n) \times \R)$.

\begin{remark}\label{remRCGITF1}
The space $(\pi_g,\taug)(\Scr{X}_f)$ is also known as {\rm the relative conormal space} of the function
$f$.
\end{remark}

\medskip
Let $f$ be as in section \ref{sectionGMTF}.

\bigskip
Since $\dim \Scr{X}_f \cap (\R^n \times \R \times \GR (n-1,n)) =n$, defining $\Scr{X}_f^\infty$ and
$\Scr{X}_t^\infty$ respectively, for a regular value $t$, as
\begin{center}
$\Scr{X}_f^\infty : = \Scr{X}_f \cap (\FH_\R^\infty \times \GR (n-1,n) \times \R)$ and $\Scr{X}_t^\infty =
\Scr{X}_f \cap (\FH_\R^\infty \times \GR (n-1,n) \times \{t\})$,
\end{center}
we deduce  $\dim \Scr{X}_f^\infty  \leqslant n-1$ and
$\dim \Scr{X}_t^\infty  \leqslant n-1$.

\medskip
Let $X_f = (\pi_f,\tf) (\Scr{X}_f)$, then $X_f$ is definable and $\dim X_f = n$. Note that $X_f$ is the
projective closure of the graph of the function $f$.

\medskip
Let $X_f^\infty = (\pi_f,\tf) (\Scr{X}_f^\infty) = X_f \cap(\FH_\R^\infty \times \R)$, thus $X_f^\infty$ is
definable and $\dim X_f^\infty \leqslant n-1$.

If $(\lambda,t) \in X_f^\infty$, let $\Om_{\lambda,t} \subset \GR (n-1,n)$ be
$\tauf((\pi_f,\tf)^{-1}(\lambda,t))$.

\begin{remark}\label{remRCGITF2}
If the subset of limits of tangent hyperplanes $\Om_{\lambda,t}$ is finite, from section
\ref{sectionGMTF}, we deduce $\# \Om_{\lambda,t} \leqslant 2 N_f$.
\end{remark}

We finally define $X_t^\infty  = \pi_f(\Scr{X}_t^\infty) \subset\FH_\R^\infty$, which is definable and of
dimension at most $n-1$. It is important to note that it may strictly contain $\FH_\R^\infty \cap
\clos(F_t)$, where $\clos (F_t) \subset \FP_\R^n$ is the projective closure of the level $F_t$.

\medskip
As a consequence of these definitions we get the following

\begin{corollary}\label{corCGI2}
(i) There exist at most finitely many values $t \in \R$ such that
$\Scr{X}_t^\infty$ is of dimension exactly $n-1$.

(ii) There exist at most finitely many points $(\lambda,t) \in X_f^\infty$ such that
$\Om_{\lambda,t}$ is of dimension $n-1$.
\end{corollary}
\begin{myproof}
Since $(X_t^\infty)_{t \in \R}$ is a definable family of subsets
whose union is  $X_f^\infty$ we then get the first point.

The second point is true for exactly the same reason for the definable family 
$(\Om_{\lambda,t})_{\{(\lambda,t)\}}$.
\end{myproof}

First, if $\Om_{\lambda,t}$ is of dimension $n-1$ so is $\Scr{X}_t^\infty$. Once more what can be said
about these values at which $\Scr{X}_t^\infty$ is of dimension exactly $n-1$ or at which
$\Om_{\lambda,t}$ is of dimension $n-1$?
%
%
%
%
%
%
%
%
%
%
%
%
%
\section{Asymptotic critical values and bifurcation values}\label{sectionACVBV}

In this section we will deal with definable functions as well as with complex polynomials. In the complex
domain we will understand differentiability in the complex meaning.

\medskip
Let us begin with the following
\begin{theorem}[\cite{Th},\cite{Ve},\cite{Har},...]\label{thmACVBV1}
Let $f : \K^n \mapsto \K$ be either a $C^l$ definable function or a complex polynomial. There exists a
smallest finite subset $B(f)$ of $\K$, called the set of bifurcation values of the function $g$ such that
for each $c \notin B(f)$ there exists an open neighbourhood $D$ of $c$ that does not meet with $B(f)$ and
such that $f_{\mid D}$ induces a $C^{l-1}$ trivial fibration over $D$, that is $f_{\mid D}^{-1}(D)$ is
$C^{l-1}$-diffeomorphic to $D \times f^{-1}(c)$.
\end{theorem}

Obviously the critical values are bifurcations values (if you do not weaken the trivialisation to be only
continuous as with the real function $f(t) = t^3$). But unfortunately, there may also exist regular values
through which the topology of the fibres is changing. In his original work using stratification theory,
Thom did not provide any means to recognise which regular value is likely to be a bifurcation value !

\medskip
For $g$ a $C^1$ function $\K^n \mapsto \K$, let $\nabla g$ be the vector field $\sum_i
\partial_{x_i}g \frac{\partial}{\partial{x_i}}$.

Let us recall what the Malgrange condition is.

\begin{definition}\label{definACVBV1}
Let $g: \K^n \mapsto \K$, be a $C^1$ function. The function $g$ satisfies the {\rm Malgrange condition at}
$c \in \K$, if for each $R\gg 1$, there exist positive constants $C$ and $\eta$ such that
\begin{center}
$x\in \{y \in \K^n : |y|> R, |f(y) - c|<\eta \} \Longrightarrow |x|\cdot |\nabla g (x)| > C$.
\end{center}
\end{definition}

\begin{definition}\label{definACVBV2}
A value $c\in \K$ is an {\rm asymptotic critical value} of the $C^1$ function $g: \K^n \mapsto \K$, if the
Malgrange condition is not satisfied at $c$, that is there exists a sequence $(x)$ of
points of $\K^n$ such that \\
(i) $|x|\rightarrow +\infty$, \\
(ii) $g(x)\rightarrow c$ and, \\
(iii) $|x| \cdot |\nabla g (x)| \rightarrow 0$ when $|x|$ goes to
$+\infty$.
\end{definition}

Let us denote by $K_\infty (g)$ the set of asymptotic critical values of the function $g$, and we define
$K(g)$ to be $K_0 (g) \cup K_\infty (g)$, the set of {\it generalised critical values}. Let us recall the
following

\begin{theorem}[\cite{Ph},\cite{HL},\cite{Pa1},\cite{Ti},\cite{LZ},\cite{dac},...]\label{thmACVBV2}
Let $f : \K^n \mapsto \K$ be either a $C^1$ definable function or a complex polynomial. \\
(1) $K(f)$ is finite. \\
(2) $B (f) \subset K (f)$.
\end{theorem}

\smallskip
Let $f$ as in Section \ref{sectionGMTF}. We also have

\begin{theorem}[\cite{DDVG3}]\label{thmACVBV3}
Let $c$ be a value. There exists a continuous definable function
germ $\theta_c : ]0,\varepsilon_c[ \mapsto ]0,+\infty[$ such that
\\
(i) there exists a constant $A >0$ such that
$\forall t \in ]0,\varepsilon_c[$, $\theta_c (t) \geqslant At$.
\\
(ii) $|x|\gg 1$ and $|f(x) - c| \ll 1$ $\Longrightarrow |x| \cdot |\nabla f (x)| \geqslant \theta_c(|f(x) -
c|)$.
\\
(iii) If $\phi$ is any definable function satisfying properties (i) and (ii), and distinct from
$\theta_c$, then there exists $0<\varepsilon<\varepsilon_c$ such
that $0<\phi (t) < \theta_c (t)$ for $t\in ]0,\varepsilon[$.
\end{theorem}

Condition (iii) also implies that, given any constant $M >1$, there exists a sequence $(x)$ going to
infinity, along which $f(x)$ tends to $c$ and such that $|x| \cdot |\nabla f (x)| \leqslant M
\theta_c(|f(x) - c|)$. Thus $c$ is an asymptotic critical value if and only if $\theta_c (t) \rightarrow 0$
when $t \rightarrow 0$.

Once an orthonormal system of coordinates is given, which we assume,
the gradient vector field $\nf$ splits into two orthogonal components, namely
its radial part $\dr f$ and its spherical part $\np f$, that is
\begin{center}
\vspace{4pt}
for $x \neq 0$ $|x| \dr f (x) = \langle \nf ,x \rangle$ and $\np f =\nf - \dr f $.
\vspace{4pt}
\end{center}

\medskip
Given $a<b$ two values, let
$\upsilon_{a,b} : ]0,+\infty[ \mapsto [0,+\infty[$ be the definable function defined as
\begin{center}
$\upsilon_{a,b} (r) = \max \left\{\displaystyle{\frac{|\dr
f|}{|\nabla f|}} : x \in  f^{-1}([a,b]) \cap \Sr_r^{n-1}\right \}$.
\end{center}

Let $\upsilon_c$ be defined as $\liminf_{\varepsilon \rightarrow 0}
\upsilon_{c - \varepsilon,c + \varepsilon}$. It is again a definable function.

\bigskip
Let us consider values $a,b$ such that $a<b$.
Let $x_0$ be a point in $f^{-1}(a)$ and let $\gamma_{x_0}$ be the trajectory of the vector
field $\frac{\nf}{|\nf|^2}$ through $x_0$ such that $\gamma_{x_0} (a) = x_0$. Thus we find

\medskip
$
\begin{array}{rcl}
\vspace{6pt}
\displaystyle{\left| |\gamma_{x_0}(t)| - |x_0| \right| = \left | \int_{a}^t
\frac{\ud|\gamma_{x_0}(\tau)| }{\ud \tau} \ud \tau \right |}
& \leqslant & \displaystyle{\int_{a}^t \left |\frac{\ud |\gamma_{x_0}(\tau)|}{\ud \tau}
\right | \ud \tau}
\\
\vspace{6pt}
& \leqslant &  \displaystyle{\int_{a}^t \left | \left \langle \frac{\ud \gamma_{x_0}
}{\ud \tau} , \frac{\gamma_{x_0}(\tau)}{|\gamma_{x_0}(\tau)|} \right \rangle
\right | \ud \tau}
\\
\vspace{6pt}
& \leqslant &  \displaystyle{\int_{a}^t \frac{|\dr f (\gamma_{x_0}(\tau))|}{|
\nf (\gamma_{x_0})|^2}\ud \tau \, .}
\end{array}
$

\begin{lemma}\label{lemACVBV1}
If $[a,b] \cap K(f) = \emptyset$, then $\upsilon_{a,b} (r)
\rightarrow 0$ when $r\rightarrow +\infty$.
\end{lemma}
\begin{myproof}
There exists a positive constant $C$ such that for $R$ large
enough, $|x| \cdot |\nf| \geqslant C$ once $x \in   f^{-1}([a,b])
\setminus \B_R^n$.

Assume there is a positive constant $A$ such that $\upsilon_{a,b}
(r) > 2A$ for $r$ large enough. We assume this is
occurring along a definable path $\alpha : [-\varepsilon,0] \mapsto
f^{-1}([a,b])\setminus \B_R^n$, such that $f\circ \alpha (s) = s+b$ and $|\alpha
(s) | \rightarrow +\infty$  as $s$ goes to $0$ and verifying
$\upsilon_{a,b} (\alpha (s)) \geqslant A$.

Taking the derivative respectively to $s$
gives $\langle \nabla f (\alpha),\alpha^ \prime \rangle = 1$. Since
\begin{center}
\vspace{3pt}
$\left \langle \displaystyle{\frac{\alpha^\prime (s)}{|\alpha^\prime (s)|}
,\frac{\alpha (s)}{|\alpha (s)|}}
\right \rangle \fl 1$ as $s \fl + \infty$,
\vspace{3pt}
\end{center}
we deduce $|\dr f| \cdot |\alpha^\prime| \fl 1$, thus

\medskip
\begin{tabular}{ccccccc}
\vspace{6pt}
$|\alpha^\prime|$ & $\leqslant$ & $\displaystyle{\frac{2}{|\dr f|}}$ &
$\leqslant$ & $\displaystyle{\frac{2}{A|\nabla f|}}$ & $\leqslant$ &
$\displaystyle{\frac{2 |\alpha (s)|}{AC}}$
\\
\vspace{6pt}
$\displaystyle{\frac{|\alpha^\prime (s)|}{|\alpha (s)|}}$ &
$\leqslant $ & $\displaystyle{\frac{2}{AC}}$ & & & &
\\
\vspace{6pt}
$\ln (|\alpha (s)|)$ & $\leqslant$ &
$\displaystyle{\frac{2}{AC} \int_{-\varepsilon}^{s} } \ud t$ & $+$ &
$\ln (|\alpha (-\varepsilon) |)$  & &
\\
\vspace{6pt}
$|\alpha (s)|$ & $\leqslant$ & $D$,  & & & &
\end{tabular}

with $D >0$ independent of $s$,
which is a contradiction to $|\alpha (s)| \fl + \infty$.
\end{myproof}

\begin{lemma}\label{lemACVBV2}
If $[a,b] \cap K(f) = \{b\}$, where $b$ is a
regular value such that $\theta_b^{-1}(t)$ is integrable when $t
\rightarrow 0$, then $\upsilon_{a,b} (r) \rightarrow 0$ when
$r\rightarrow +\infty$.
\end{lemma}
\begin{myproof}
Assume there is a positive constant $A$ such that $\upsilon_{a,b} (r) >
2A$ for $r$ large enough.
We can assume this phenomenon is occurring along a definable path
$\alpha : [-\varepsilon,0[ \mapsto f^{-1}([a,b])$,
such that $f\circ (\alpha (s)) = s+b$ and $|\alpha (s) | \rightarrow
+\infty$  as $s$ goes to $0$ and verifying $\upsilon_{a,b} (\alpha (s))
\geqslant A$.
We also assume that $\varepsilon$ is such that for each $s > -\varepsilon$, $|\alpha (s)|$
satisfies the point (ii) of Theorem \ref{thmACVBV3} at $b$.

Taking the derivative respectively to $s$
gives $\langle \nabla f (\alpha),\alpha^ \prime \rangle = 1$. Since
\begin{center}
\vspace{4pt}
$\left \langle \displaystyle{\frac{\alpha^\prime (s)}{|\alpha^\prime (s)|}
,\frac{\alpha (s)}{|\alpha (s)|}}
\right \rangle \fl 1$ as $s \fl + \infty$,
\vspace{4pt}
\end{center}
we deduce $|\dr f| \cdot |\alpha^\prime| \fl 1$, thus

\medskip
\begin{tabular}{ccccl}
\vspace{6pt}
$|\alpha^\prime|$ & $\leqslant$ & $\displaystyle{\frac{2}{|\dr f|}}$ &
$\leqslant$ & $\displaystyle{\frac{2}{A|\nabla f|}} \; \leqslant \;
\displaystyle{\frac{B |\alpha (s)|}{|\theta_b (-s)|}} \;$ for a constant
$B>0$
\\
\vspace{6pt} &  & $\displaystyle{\frac{|\alpha^\prime (s)|}{|\alpha (s)|}}$ & $\leqslant $ &
$\displaystyle{\frac{B}{|\theta_b(-s)|}}$
\\
\vspace{6pt}
 &  & $\ln (|\alpha (s) |)$ & $\leqslant$ & $B \displaystyle{\int_{-\varepsilon}^s
\frac{1}{|\theta_b(-t)|}} \ud t$  + const.
\\
\vspace{6pt}
 & & $|\alpha (s)|$ & $\leqslant$ & $C$, with $C >0$ independent of $s$,
\end{tabular}

which is a contradiction to $|\alpha (s)| \fl + \infty$.
\end{myproof}

So from this we recover the following
\begin{proposition}[\cite{DDVG3}]\label{propACVBV1}
Let $c$ be a regular value. If $\theta_c^{-1}$ is integrable nearby $0$, then $c \notin B(f)$ and
the trivialisation is realised by the local flow of $|\nf|^{-2}\cdot \nf$.
\end{proposition}

\medskip
When $f$ is a real or complex polynomial (or more generally definable in a polynomially bounded o-minimal
structure), given any value $c$ there exist a smallest real number $\rho_c$ (belonging to the field of
exponents of the o-minimal structure) and a positive constant $L_c$ such that
\begin{center}
\vspace{4pt}
if $|x| \gg 1$ and $|f(x) - c| \ll 1$ then $|x| \cdot |\nabla f(x)| \geqslant L_c|f(x) -c|^{\rho_c}$,
\vspace{4pt}
\end{center}
that is
\begin{center}
$\theta_c (t) = L_c t^{\rho_c}$ and thus $\rho_c \leqslant 1$, \vspace{3pt}
\vspace{4pt}
\end{center}
from \cite{DDVG2} and \cite{DDVG3}. Thus $c$ is an asymptotic critical value if and only if $\rho_c >0$.
So, the integrability of $\theta_c$ is equivalent to $\rho_c <1$. Note that also that $\rho_c <1$ if and
only if $\upsilon_c (r) \fl 0$ as $r \fl + \infty$. If $c$ is a regular bifurcation value then $\rho_c
= 1$.

Let us mention another result of the same kind
\begin{theorem}[\cite{LZ}]
Let $c$ be a regular value such that $\upsilon_c <1$. Then $c \notin B(f)$.
\end{theorem}

See also \cite{NZ} for the complex polynomial version of this result.

\medskip
The first remark is that the trivialisation is provided by a vector field tangent to the spheres that are
transverse to the levels of $f$ nearby $c$ near infinity. Loi and Zaharia gave a more general version: They
require transversality to the levels of a proper definable positive $C^1$ submersion $\mu$ on $\Rn
\setminus\{0\}$ \cite{LZ}. They proved that $B(f) \setminus K_0 (f) \subset S_\mu (f) \subset K_\infty
(f)$, where $S_\mu (f)$ is defined as
\begin{center}
\vspace{4pt}
$ S_\mu (f) :=\left\{b \in \R: \exists (x) : |x| \fl +\infty, f(x) \fl b, \left\langle\nu_f (x),
\displaystyle{\frac{\nabla\mu (x)}{|\nabla\mu (x)|}} \right\rangle=\pm 1 \right\}$.
\vspace{4pt}
\end{center}

To finish this section let us relate with Section \ref{sectionRCGITF}.
\begin{proposition}\label{propACVBV2}
Let $c$ be a regular value of $f$ such that there is $\lambda \in \FH_\R^\infty$ such that $\Om_{\lambda,c}$ is
of dimension $n-1$. Then $c \in K_\infty (f)$.
\end{proposition}
\begin{myproof}
Let $\xi_\lambda$ be a unit vector collinear to the line direction $\lambda$.
Thus there exists an open subset $\Omega \subset \Srn$ such that for any $u \in \Omega$, there
exists a connected component $\Gamma_u$ of the oriented polar curve $\Psi_f^{-1} (\{u\} \times\R)$ such that
$\Gamma_u$ is unbounded, and $f_{\mid \Gamma_u} (x) \fl c$ as $\Gamma_u \ni x \fl \infty$ such that
$x /|x| \fl \pm \xi_\lambda$. Thus there is a unit vector $u$ such that $\langle u, \xi_\lambda\rangle \neq 0$. Thus
for $r$ large enough, $2\upsilon_c (r) \geqslant |\langle u, \xi_\lambda\rangle|$.
\end{myproof}

The converse of this result is not true in the real context \cite [example 5.3]{DDVG2}. In Section
\ref{sectionRPCP} we will see that in the complex polynomial case, this point deserves to be discussed.
%
%
%
%
%
%
%
%
%
%
%
%
%
%
\section{Triviality at infinity of tame functions with Strongly Isolated Singularity
at Infinity}\label{sectionTITFSISI}

We still assume that $f$ is as in Section \ref{sectionGMTF}. For our purpose here, we
do the extra assumption that $K_0 (f) = \emptyset$ and $K(f)
\subset \{c\}$, since we will only deal with asymptotic critical values that are regular values.

In this section we give a sufficient condition to trivialise the function $f$ in a neighbourhood of the
regular value $c$ that will be expressed in terms of the total absolute curvature.

\medskip
We recall from Section \ref{sectionRCGITF} that $\FH_\R^\infty = \FP_\R^n \setminus \Rn$, so we also
consider any subset of $\G_\R (1,n)$ as a subset of $\FH_\R^\infty$ if needed.

Let us introduce the sufficient condition we just mentioned above, and that we have called the {\bf SISI} 
condition.
\begin{definition}\label{definTITFSISI1}
Let $f$ as above. Let $c$ be a regular value taken by $f$. The function $f$ is said to have {\rm strongly
isolated singularities at infinity} at $c$ if the following condition is satisfied:
\\
There exists a finite subset $\Lambda_c \subset \GR (1,n)$ such that for each line direction $\lambda \in
X_c^\infty \setminus \Lambda_c$, for each hyperplane direction $\bh \in \Om_{(\lambda,c)}$, the line
direction $\lambda$ is contained in $\GR(1,\bh)$ the Grassmann space of line directions of the hyperplane
$\bh$.
\end{definition}

The main result of the paper is the following

\begin{theorem}\label{thmTITFSISI1}
Let $f:\Rn\mapsto\R$ be a $C^l$ definable function with $l\geqslant 2$. Assume that
the function $f$ satisfies condition {\bf SISI} at $c$.
\\
If the function  $t \mapsto |K|(t)$ is continuous at $c$, the function $f$ is trivialised over a
neighbourhood of $c$ by means of the flow of a $C^{l-1}$ definable vector field. So $c \notin B(f)$.
\end{theorem}

The rest of the section is devoted to the proof of this result.

\medskip
To prove Theorem \ref{thmTITFSISI1} we first need the following
\begin{lemma}\label{lemTITFSISI1}
Under the hypotheses of Theorem \ref{thmTITFSISI1},  $\tauf(\Scr{X}_c^\infty)$ is of dimension at most
$n-2$.
\end{lemma}
\begin{myproof}
From Proposition \ref{propGMTF1} and the continuity of the total absolute curvature function at $c$ we
deduce that $\cup_{\lambda \in X_c^\infty} \Om_{(\lambda,c)}$ is of dimension at most $n-2$.
\end{myproof}

\begin{myproof} 
Let $X_c^{\infty,+}$ be the lift of $X_c^\infty$ onto $\Srn$. It is a closed definable subset of dimension 
at most $n-2$, let $\Lambda_c^+$ be the lift of $\Lambda_c$.

\medskip
After a direct orthonormal change of coordinates if necessary, we can assume that the intersections of each
coordinate axis with the unit sphere does meet $X_c^{\infty,+}$ and that any $u \in \Lambda_c^+$ does not lie
in any coordinates hyperplane.

\medskip
Let $\Delta :=\{(\delta_1,\ldots,\delta_n): \forall i, \delta_i >0\}$. 
Embedding $\Delta$ in ${\rm Gl}_n(\R)$ as diagonal matrices
makes $\Delta$ a smooth semi-algebraic subgroup with a smooth semi-algebraic action over $\Rn$.
Note that $\Delta$ is diffeomorphic to $\R^n$.

For any $A \in \Delta$, let us define the following semi-algebraic function
\begin{center}
\vspace{3pt}
$g_A(x) = \langle A\cdot x,x\rangle^\frac{1}{2}$.
\vspace{3pt}
\end{center}

Note that $g_A$ is a smooth proper submersion outside the origin since such an $A$ is positive definite.

We will recycle here the method used in \cite{NZ} and \cite{TZ} but with "spheres" given 
by the levels of a function $g_A$ for an appropriate $A$.

\medskip
For each $u \in \Lambda_c^+$, let $V(u)$ be the closed definable subset defined as
\begin{center}
\vspace{4pt}
$V(u) :=\clos\{\nu \in \Srn: \nu = \lim \nu_f (x)$ with $|x| \fl +\infty$, $\displaystyle{\frac{x}{|x|}\fl u}
\hfill$

$\hfill$ and $f(x)\fl c\}$.
\vspace{4pt}
\end{center}
Obviously each $V(u)$, $u \in \Lambda_c^+$, is of dimension at most $n-2$.

\medskip
Let $\Delta_0 :=\{A\in \Delta: i\neq j \Longrightarrow \delta_i \neq \delta_j\}$. It is open, semi-algebraic
and dense in $\Delta$.
Since in the new coordinate system any $u \in X_c^{\infty,+}$ has at least two non zero coordinates, 
no such vector $u$ can be an eigenvector of $A \in \Delta_0$.

As a corollary of this fact we get
\begin{lemma}\label{lemTITFSISI2}
For each $A \in \Delta_0$, there exists $\alpha,\beta \in ]0,1[$ such that there are $R>0$ and $\varepsilon>0$ 
such that for each $x \in f^{-1}(]c-\varepsilon,c+\varepsilon[) \setminus \clos (\B_R^n)$
\begin{center}
\vspace{4pt}
$\alpha <\displaystyle{\left\langle \frac{A\cdot x}{|A\cdot x|},\frac{x}{|x|}\right\rangle}<\beta$.
\vspace{4pt}
\end{center}
\end{lemma}
\begin{myproof}
Since $X_c^{\infty,+}$ is compact and by definition of $\Delta_0$, there exists positive $\alpha_0$ 
and $\beta_0$ such that for any $u\in X_c^{\infty,+}$,
\begin{center}
\vspace{4pt}
$\alpha_0 <\displaystyle{\left\langle \frac{A\cdot u}{|A\cdot u|},u\right\rangle}<\beta_0<1$.
\vspace{4pt}
\end{center}
Writing the definition of $X_c^{\infty,+}$ provides the desired statement.
\end{myproof}

Given $u\in\Lambda_c^+$, each coordinates is non zero, so the subset $\Delta_0 (u):= \{A\cdot u: A \in 
\Delta_0\}$ is a semi-algebraic open subset of $\Rn\setminus\{0\}$. So its image under the radial 
projection is open in $\Srn$. Let $\rho : \Rn\setminus\{0\} \mapsto \Srn$ be the radial projection. 
It is a smooth semi-algebraic map thus 
$U_{c\in\Lambda_c^+} \rho^{-1} (V(u))$ is a definable positive cone of dimension at most $n-1$.
This means there exists $\Delta_1$ an open dense definable subset of $\Delta$ such that for any 
$u\in\Lambda_c^+$ and for any $A\in\Delta_1$, $A \cdot u \notin V(u)$.

Given $u\in \Srn$ an $\eta >0$, the positive conical neighbourhood of the oriented semi-line 
$\R^+ u$ of radius $\eta$ is the following
\begin{center}
\vspace{4pt}
$C^+(u;\eta):=\{x \in \Rn\setminus\{0\}: |u - \rho (x)|<\eta\} \cup \{0\}$.
\end{center}

\begin{lemma}\label{lemTITFSISI3}
Let $A \in \Delta_1$. There exists $\gamma,\delta \in ]0,1[$ such that there are $R>0$ and $\varepsilon>0$, 
$\eta>0$, such that for each $x \in (f^{-1}(]c-\varepsilon,c+\varepsilon[) \cap C^+(u;\eta)) \setminus 
\clos (\B_R^n)$,
\begin{center}
\vspace{4pt}
$-1<-\gamma <\displaystyle{\left\langle \frac{A\cdot x}{|A\cdot x|},\frac{\nf (x)}{|\nf (x)|}\right\rangle}<
\delta<1$.
\vspace{4pt}
\end{center}
\end{lemma}
\begin{myproof}
From the definition of $\Delta_1$ we deduce, there exists $\gamma_0,\delta_0 \in ]0,1[$ such that for 
each $u\in \Lambda_c^+$, each $\nu \in V(u)$ and each $v \in \Srn$ such that $|v-u|<\eta_0$ for some 
positive $\eta_0$, we find 
\begin{center}
\vspace{4pt}
$-1<-\gamma_0 <\displaystyle{\left\langle \frac{A\cdot v}{|A\cdot v|},\nu\right\rangle}<\delta_0<1$.
\vspace{4pt}
\end{center}
Returning to the definitions of $X_c^{\infty,+}$ and $V(u)$ gives us the desired uniform version.
\end{myproof}

We recall that $\nu_{g_A}(x) = \displaystyle{\frac{\nabla g_A (x)}{|\nabla g_A (x)|} = 
\frac{A\cdot x}{|A\cdot x|}}$.
As a corollary of Lemma \ref{lemTITFSISI2} and of Lemma \ref{lemTITFSISI3} we obtain
\begin{proposition}\label{propTITFSISI1}
For each $A \in \Delta_1$, there exists $\alpha,\beta,\gamma \in ]0,1[$ such that there 
exist $R >0$, $\varepsilon>0$ and $\eta>0$ such that for each $x \in f^{-1}(]c-\varepsilon,
c+\varepsilon[) \setminus \clos (\B_R^n)$

(i) $\alpha |x|<\langle \nu_{g_A}(x),x\rangle<\beta|x|$,

(ii) for each $u\in \Lambda_c^+$ such that $x \in C^+(u;\eta)$, $|\langle \nu_{g_A}(x),\nu_f (x)\rangle|<
\gamma<1$.
\end{proposition}

Let $A$ be given as in Proposition \ref{propTITFSISI1}.
Assume that $R$ and $\varepsilon$ are given.
Let us define the following vector field
\begin{center}
\vspace{4pt}
$\omega_A  (x) := \nu_f(x) - \langle \nu_f (x), \nu_{g_A} (x)\rangle\nu_{g_A}(x)$.
\vspace{4pt}
\end{center}
For $x  \in f^{-1}(]c-\varepsilon,c+\varepsilon[) \setminus \clos (\B_R^n)$ this vector field is 
non vanishing, since $|\omega_A (x)| \geqslant \sqrt{1-\gamma^2}$, and is tangent to the levels of 
$g_A$ (which are compact) and is transverse to the levels of $f$. 

\medskip
Let $\delta$ be the biggest eigenvalue of $A$.
Now we define the vector field that will realise the trivialisation around $c$: \\

For $x$ such that $g_A(x) \leqslant \delta R$,  let $\xi (x) := \nu_f(x)$, \\

for $x$ such that $g_A(x) \geqslant 2\delta R$, let
$\xi (x) := \displaystyle{\frac{\omega_A (x)}{|\omega_A (x)|}}$, \\

for $g_A^{-1}(x) \in[\delta R,2\delta R]$, let

\begin{center}
$\xi (x) := \displaystyle{\kappa(g_A(x))\nu_f(x)+ [1 - \kappa(g_A(x))]
\frac{\omega_A(x)}{|\omega_A (x)|}}$,
\vspace{4pt}
\end{center}
where $\kappa: [\delta R,2\delta R] \mapsto [0,1]$ is a $C^{l}$ definable function that is strictly 
decreasing and such that $\kappa(\delta R) = 1$ and $\kappa(2\delta R)=0$, and is also 
$l$-flat at $\delta R$ and at $2\delta R$.

\medskip
Restricting $\xi$ to $f^{-1}(]c-\varepsilon,c+\varepsilon[)$, we observe that $\xi (x)$ does not 
vanish and so the trivialisation of $f$ in a neighbourhood of $c$ is provided by the 
flow of $\xi$ as in \cite{NZ} and \cite{LZ}.
\end{myproof}

\section{Real polynomials versus complex polynomials}\label{sectionRPCP}

Let $\K$ be either $\R$ or $\C$. Let $f : \K^n \mapsto \K$ be a polynomial of degree $d \geqslant 2$. Let
$\nf$ be the polynomial vector field $\sum_i \partial_{x_i} f \partial_{x_i}$.

\medskip
If $|\nf (x)| \fl 0$ along a sequence $x$, such that $|x| \fl +\infty$ and $x /|x| \fl \bu \in
\FP_\K^{n-1}$, we deduce that for each $i = 1,\ldots,n$, $\partial_{x_i} f_d (\bu) = 0$, when $f =f_d +
f_{d-1} + \ldots + f(0)$ is written as the sum of its homogeneous components. Assume moreover that along
this sequence the Malgrange condition fails at $c$ a regular value, that is $|x|\cdot|\nf (x)| \fl 0$ and
$f(x) \fl c$.

Assume, after a rotation, that $\bu = (0,\ldots,0,1)$, then
writing $y_i = x_i /x_n$ and $y_0 = 1/x_n$, we deduce that along the sequence $x$
\begin{center}
\vspace{6pt}
$|(\partial_{y_1} \ft^{(n)},\ldots,\partial_{y_{n-1}} \ft^{(n)},d\ft^{(n)}-y_0\partial_{y_0} \ft^{(n)})|
\ll |y_0|^d$,
\end{center}
with
\begin{center}
$\ft^{(n)} (y_0,\ldots,y_{n-1}) = y_0^d f(y_1/y_0,\ldots,y_{n-1}/y_0,1/y_0)$.
\vspace{6pt}
\end{center}
We deduce that for $i = 1,\ldots,n-1$, $\partial_{y_i} f_d (\bu) = 0$, and so $\partial_{x_n} f_d (\bu) =0$.
So we also find that $|\partial_{y_0} \ft^{(n)}| \leqslant const \cdot |y_0|^{d-1}$,
thus $f_{d-1} (0,\ldots,0,1) = f_{d-1}(\bu) = 0$.

\medskip
These elementary computations mean that the set of points $\bu$ at infinity nearby which Malgrange
condition fails at $c$ is very specific, namely these points are roots of $f_{d-1}$ and of $\nabla f_d$.

\medskip
Let $\lambda \in \G_\K (1,n)$ be a line direction and let $\bh_\lambda \in \G_\K(n-1,n)$ be the hyperplane
direction orthogonal to $\lambda$. We denote by ${\bf P}_\lambda (f)$ the subset $\{x \notin \crit (f): T_x
f = \bh_\lambda\}$, that is the polar variety of the function $f$ in the line direction $\lambda$. Note that
it is a semi-algebraic subset of $\K^n$.

\begin{lemma}\label{lemRPCP1}
Let $c$ be a regular value of $f$. Assume there exists a line direction $\lambda$ such that there exists a
sequence $(x)$ in ${\bf P}_\lambda (f)$, $|x|\fl +\infty$, $x/|x| \fl \bu$ and $f(x) \fl c$. If
$f_{d-1}(\bu) \neq 0$, then  $\nabla f_d (\bu) \neq 0$ and so $\lambda = \R \nu_{f_d} (\bu)$.
\end{lemma}
\begin{myproof}
We can assume that this phenomenon occurs along a semi-algebraic path, say $\Gamma \subset {\bf P}_\lambda
(f)$ and is parametrised as $]0,\varepsilon[ \ni r \mapsto x(r)$ such that $|x(r)| \fl + \infty$ as $r \fl 0$
and $r |x(r)| \fl 1$.
\\
Assume that $f_{d-1} (\bu) = \alpha_0 \neq 0$. Thus there exists a positive constant $M$ 
such that for $r$ small enough
\begin{center}
\vspace{4pt}
$|x(r)| \cdot |\nf (x(r))| \geqslant M$, that is $|\nf (x(r))| \geqslant M \cdot  r $.
\vspace{4pt}
\end{center}
Since along the path $|\dr f| \sim |x| \cdot |f-c|$, we deduce that $|\np f| \gg |x| \cdot |f-c|$. 
From Lemma \ref{lemACVBV1} we get $\langle \bu,\lambda \rangle = 0$.

\medskip
After an orthonormal change of coordinates we  assume that $\bu = (0,\ldots,0,1)$. \\
Writing $\lambda =(\lambda_1,\ldots,\lambda_n)$ and since $\langle \bu, \lambda\rangle = 0$,
after a rotation in $\lambda_1,\ldots,\lambda_{n-1}$ we actually get $\lambda_1 = \pm 1$ and
$\lambda_2 = \ldots = \lambda_n = 0$. So we get  $(\partial_{x_i} f)_{\mid \Gamma} = 0$, $i=2,\ldots,n$.

\medskip
Let $y_i := x_i / x_n$, for $i=1,\ldots,n-1$, and $y_0 := 1/x_n$. Thus we get
\begin{center}
\vspace{3pt}
$\ft^{(n)} (y_0,\ldots,y_{n-1}) := y_0^d f(y_1/y_0,\ldots,y_{n-1}/y_0,1/y_0)$.
\vspace{3pt}
\end{center}
By abuse of notation let us define
\begin{center}
\vspace{3pt}
$\partial_{y_n} \ft^{(n)} (y_0,\ldots,y_{n-1}) := y_0^{d-1}
\partial_{x_n} f(y_1/y_0,\ldots,y_{n-1}/y_0,1/y_0)$.
\vspace{3pt}
\end{center}
Thus we get $d\cdot\ft^{(n)} =  \partial_{y_n} \ft^{(n)}  + \sum_{i=0}^{n-1} y_i \partial_{y_i} \ft^{(n)}$, and so
\begin{center}
\vspace{3pt}
$|y_0|^{d-1} \cdot|\nf (x)| = |(\partial_{y_1} \ft^{(n)}(y),\ldots,\partial_{y_n} \ft^{(n)} (y))|$.
\vspace{3pt}
\end{center}

Let $y(r) = (y_0(r),\ldots,y_{n-1}(r))$ be the path in the new coordinates.
Thus we get $y_0 \simeq r$ and $r^{-1}y_i (r) \fl 0$ as $r$ goes to $0$. Thus we can assume $y_0 =r$.

Thus along $\Gamma$ we know  that $r^{-d} \cdot|\ft^{(n)} (y(r)) - cr^d| \fl 0$ as $r$ goes to $0$, and
\begin{center}
\vspace{3pt}
$|(\partial_{y_1} \ft^{(n)}(y),\ldots, \partial_{y_n} \ft^{(n)} (y))| \geqslant Mr^d$, that is
$|\partial_{y_1} \ft^{(n)}(y)| \geqslant Mr^d$.
\vspace{3pt}
\end{center}
Along $\Gamma$ we deduce that
\begin{center}
\vspace{3pt}
$d \cdot \ft^{(n)} (y(r)) = r \partial_{y_0} \ft^{(n)} (y(r)) + y_1
\partial_{y_1} \ft^{(n)} (y(r)) = \hfill$ \\
 $\hfill rf_{d-1} (y(r)) + y_1 (r) \partial_{y_1} f_d(y(r)) + o(r)$.
\vspace{3pt}
\end{center}
Taking the derivative of $\ft^{(n)}(y(r))$ in $r$ provides
\begin{center}
\vspace{3pt}
$\partial_{y_0} \ft^{(n)}  + y_1^\prime \partial_{y_1} \ft^{(n)} \simeq dcr^{d-1}$.
\vspace{3pt}
\end{center}
So if $\partial_{y_1} \ft^{(n)} \simeq \alpha r^a$ and $y_1 (r) \simeq \beta r^b$
and since $d\geqslant 2$ we deduce that $a+b = 1$ and $\alpha_0 + b\alpha\beta = 0$.
But we also deduce that $\alpha_0 + \beta \alpha = 0$ so $b=1$ and thus $a = 0$.
Which implies that $\partial_{y_1} f_d (\bu) = \alpha \neq 0$ and so the claim is proved.
\end{myproof}

\smallskip
Let us recall what is happening in the complex case. Parusi\'nski defined in \cite{Pa1} the notion of
complex polynomial with isolated singularities at infinity. By this he means that the subset
\begin{center}
\vspace{3pt}
$A:=\{ \lambda \in \FH_\C^\infty : \partial_{x_1} f_d (\lambda)= \ldots = \partial_{x_n} f_d (\lambda) =
f_{d-1} (\lambda)= 0\}$,
\vspace{3pt}
\end{center}
is finite. For each $t \in \C$, let ${\bf F}_t \subset  \FP_\C^n$ be the projective closure of the level
$f^{-1} (t)$. So for each regular value $t$, the projective hypersurface ${\bf F}_t$ has only isolated
singularities. Let $\mu({\bf F}_t)$ be the sum of the Milnor numbers of the isolated singularities of ${\bf
F}_t$.

The notion of relative conormal bundle and relative conormal space, defined in Section \ref{sectionRCGITF},
still make sense for a complex polynomial when dealing with the complex analogs of the notions used in the
real setting. Note then, once $t \notin K_\infty (f)$, we find $X_t = {\bf F}_t$.

\medskip
For a given line direction $v$, let $\clos ( {\bf P}_\lambda (f))$ be the projective closure of ${\bf
P}_\lambda(f) \subset  \C^n \setminus \crit (f) \subset \FP_\C^n$. Then

\begin{theorem}[\cite{Pa1,Ti}]\label{thmRPCP1} Let $c$ be a regular value of the complex polynomial $f$
with isolated singularities at infinity. The following statements are equivalent:

(1) $c \notin K_\infty (f)$.

(2) $c \notin B(f)$.

(3) The total Milnor number function $t \mapsto \mu ({\bf F}_t)$ is locally constant in a neighbourhood of
$c$.

(4) The Euler Characteristic function  $t \mapsto \chi ({\bf F}_t)$ is locally constant in a neighbourhood
of $c$.

(5) For $\lambda$ in a Zariski open set of $\G_\C(1,n)$, $\clos ({\bf P}_\lambda (f)) \cap X_c^\infty =
\emptyset$.

(6) The dimension of $\Scr{X}_c^\infty$ is at most $n-2$.
 \end{theorem}

\medskip
As far as the author knows, there is no real polynomial version of such a statement. One of the reason for
this is that the local constancy of simple invariants such as the Euler Characteristic or some Milnor
numbers is not usually a sufficient condition to ensure the equisingularity of a family.

\medskip
Let us mention the connection between Theorem \ref{thmRPCP1} and the complex version of Theorem
\ref{thmACVBV1}, which holds true with $\theta_c (t) = K_c t^{\rho_c}$ for a rational number
$\rho_c\leqslant 1$ and $K_c >0$.

\begin{proposition}[\cite{DDVG2}]\label{propRPCP1}
Let $f$ be a complex polynomial with isolated singularities at infinity. Let $c$ be a regular value. The
exponent $\rho_c$ is equal to $1$ if and only if $c \in B(f)$.
\end{proposition}

The real version of this result is not true, as given by $f(x,y) = - y(2x^2y^2 -9xy+12)$ (see \cite{DDVG1,DDVG2}
for more on this example).

\medskip
For real polynomial functions we have the following

\begin{proposition}\label{propRPCP2}
Let $f$ be a real polynomial function on $\Rn$ of degree $d$. Assume the following subset
\begin{center}
\vspace{3pt}
$A := \{ \lambda \in \FH_\R^\infty : \partial_{x_1} f_d (\lambda)= \ldots = \partial_{x_n} f_d (\lambda) =
f_{d-1} (\lambda) =0\}$
\vspace{3pt}
\end{center}
is finite. Given any sequence $x$ such that $|x| \fl +\infty$, $f(x) \fl c \notin K_0(f)$ and $x /|x| \fl
\lambda \notin A$, any limit of tangent hyperplane direction $\bh = \lim T_x F$ contains the line direction
$\lambda$, so the function $f$ satisfies condition {\bf SISI} at $c$.
\end{proposition}
\begin{myproof}
Let $r \mapsto r u(r)$ be a semi-algebraic path such that $|u(r)| = 1$ and $f(r u(r)) \fl c$ as $r \fl +
\infty$. Assume that $u(r) \fl \bu \in \Srn$ as $r \fl +\infty$ such that $f_{d-1} (\bu) = \alpha_0\neq 0$.
We must either have $f_d (\bu) =0$ or $\langle \nu_{f_d} (\bu) ,\bu \rangle = 0$. After a rotation we
assume that $\bu =(0,\ldots,0,1)$. Writing $u(r) = (u_1(r),u_2(r),\ldots,1-u_n(r))$, we find that $u_i
(r)\sim r^{-e_i}$ for rational positive numbers $e_i$, $i=1,\ldots,n$. Let $e = \min \{e_i\}$. Then we find
that $e_n = 2e$. For each $i = 1,\ldots n$, we get $\partial_{x_i} f_d (u(r)) \sim r^{-d_i}$, for some
positive rational numbers $d_i$. Since $f(ru(r)) \fl c$, we must have $|rf_d(u(r)) +f_{d-1}(u(r))| \leqslant
const \cdot r^{-1}$. So $2e \leqslant 1$ and $f_d (u(r)) \simeq -\alpha_0 r^{-1}$. We deduce that
\begin{center}
\vspace{3pt}
$\langle \nabla f_d (u(r)),u(r)\rangle \simeq -d \alpha_0 r^{-1}$ and $\langle \nabla f_d (u(r)),u^\prime(r)
\rangle \simeq d\alpha_0 r^{-2}$.
\vspace{3pt}
\end{center}
If $d_n \geqslant 1$, then
\begin{center}
\vspace{3pt}
 $\langle \nabla f_d (u(r)),u^\prime(r)
\rangle \simeq u_2^\prime (r) \partial_{x_2} f_d(u(r)) + \ldots + u_n^\prime (r)\partial_{x_n} f_d (u(r))
\simeq d\alpha_0 r^{-2}$.
\vspace{3pt}
\end{center}
There exists $i \in \{2,\ldots,n\}$ such that $e_i + d_i \leqslant 1$. So
\begin{center}
\vspace{3pt}
$\partial_{x_i} f_d (u(r))+r^{-1}\partial_{x_i} f_{d-1} (u(r)) \sim r^{-d_i}$
\vspace{3pt}
\end{center}
thus $\nf (ru(r)) \simeq \nabla f_d (ru(r))$ and moreover $\langle \nabla f_d/|\nabla f_d|,
\bu \rangle \fl 0$.

\noindent 
\smallskip 
If $d_n < 1$, there exists $i \in \{2,\ldots,n\}$ such that $e_i + d_i \leqslant d_n$
and thus again we deduce $\nf (ru(r)) \simeq \nabla f_d (ru(r))$ and $\langle \nabla f_d/|\nabla f_d| , \bu
\rangle \fl 0$. This ends the proof.
\end{myproof}

So condition {\bf SISI} is a reasonable condition to work with in the frame we are given.
The finiteness of the subset $A$ has also another consequence.
As a corollary we obtain
\begin{proposition}\label{propRPCP3}
Let $f$ be as in Proposition \ref{propRPCP2}. \\
(1) If $c \notin K(f)$, then the function $t \mapsto |K|(t)$ is continuous at $c$. \\
(2) If the total curvature function $|K|$ is not continuous at a regular value $c$, then
$c \in K_\infty (f)$ and moreover its exponent at infinity is $1$.
\end{proposition}
\begin{myproof}
Point (1) comes from the proof of Lemma \ref{lemRPCP1}. \\
For point (2), the discontinuity of the total curvature function at $c$ guarantees we can apply
Lemma \ref{lemRPCP1}. So there exists $v \in \Srn$ such that along the oriented polar curve in the oriented
direction $v$, namely $\Psi_f^{-1}(\{v\}\times \R)$ is not empty and moreover  there exists a branch
$\Gamma_v$ of this polar curve such that $f_{\mid \Gamma_v} (x) \fl c$ as $|x| \fl +\infty$ along which
$|(\nf)_{\mid \Gamma_v}| (x) \sim |\dr f_{\mid \Gamma_v}| (x)$, that is $\rho_c =1$.
\end{myproof}

This partially answers, in the real polynomial case, the question about the values at which the total
absolute curvature is not continuous: When regular they can only be asymptotic critical values $c$ with
$\rho_c =1$.

As a final corollary we deduce

\begin{corollary}\label{corRPCP2}
Let $f$ be a real polynomial function on $\Rn$ of degree $d$ such that
\begin{center}
\vspace{3pt}
$A := \{ \lambda \in \FH_\R^\infty : \partial_{x_1} f_d (\lambda)= \ldots = \partial_{x_n} f_d (\lambda) =
f_{d-1} (\lambda) =0\}$
\vspace{3pt}
\end{center}
is finite. If $c$ is a regular bifurcation value of $f$, then the function total curvature of $f$ is not
continuous at $c$, or equivalently there exists a non empty open subset $\Lambda \subset \GR (1,n)$ such
that for each $\lambda \in \Lambda$, ${\bf P}_\lambda (f)$ is a non empty smooth curve and $\clos ({\bf
P}_\lambda (f)) \cap X_c^\infty \cap A \neq \emptyset$.
\end{corollary}
%
%
%
%
%
%
%
%
%
%
%
\section{Comments and remarks}\label{sectionCR}
The reader will have noticed that the conclusion of Theorem \ref{thmTITFSISI1} still holds true if we drop
the hypothesis on the continuity of the total absolute curvature function, to only requiring that $\tauf
(\Scr{X}_c^\infty)$ - the closure of the limits of the tangent spaces to the fibres as they tend to $c$ -
is a proper closed subset of $\G_\R (n-1,n)$. The proof works the same. But as noticed
by Tib\u{a}r in the real polynomial case \cite{Ti}, if $\tauf (\Scr{X}_c^\infty)$ is of dimension $n-1$ then it
is $\G_\R (n-1,n)$. In this situation finer conditions (yet unknown) will be required to ensure the
trivialisation, even in the case of {\bf SISI}.

\medskip
We have said that Theorem \ref{thmTITFSISI1} was the real counterpart of Theorem \ref{thmRPCP1}. Our result
only provides a sufficient condition expressed in terms of total absolute curvature (or in terms of polar
curves) once specified the points at infinity (that is in $\FH_\R^\infty$) at which the Gauss-Kronecker
curvature may concentrate.

\medskip
In the complex domain, to each $t \in \C$, we associate a real number $LK (t)$, which is just the total
$(2n-2)$-Lipschitz-Killing curvature of the regular part of the real algebraic subset $F_t \subset \C^n$
(see \cite{TS} for how this can be used in equisingularity problems). It is of constant sign and it is
obvious that point (vi) of Theorem \ref{thmRPCP1} is equivalent to the continuity at $c$ of the function $t
\mapsto LK (t)$. So in a sense, since it is equivalent to point (iv), Theorem \ref{thmRPCP1} can be also
considered as a Gauss-Bonnet-Chern type constancy result.

\medskip
Having isolated singularities at infinity in the complex case is almost exactly requiring which points at
infinity are likely to concentrate curvature. The finiteness of $\{\partial_{x_1} f_d= \ldots =
\partial_{x_n} f_1= f_{d-1}=0\} \cap \FH_\C^\infty$ implies that if there are generic polar curves along a non
bounded branch of which $f$ tends to a regular value $c$, this branch must tend to a point in
$\{\partial_{x_1} f_d= \ldots = \partial_{x_n} f_1= f_{d-1}=0\} \cap \FH_\C^\infty$. Thus, the knowledge of
the set of points at infinity where the generic polar curves (or the "generic" non-empty polar varieties
when the Gauss map is degenerate) are ending seems to be an interesting object to understand and to
describe when we are willing to decide whether a regular value is a bifurcation value or not, see for
instance \cite[Example 2.13]{Ti}.

\bigskip
In the real domain, the polynomial case is already delicate since point (i) and point (ii) of
Theorem \ref{thmRPCP1} are already not equivalent (see \cite{TZ}).
Moreover even among regular asymptotic critical values there are distinctions to make as suggested
in Section \ref{sectionACVBV}.

King-Zaharia-Tib\u{a}r example $f(x,y) = - y(2x^2y^2 -9xy+12)$ in the real plane,
which satisfies our hypotheses, has no bifurcation value, but $0 \in K_\infty (f)$ and
$\lim_{t \fl 0}|K(t)| = 2\pi$ while $|K(0)| = 0$. The trivialisation cannot be realised by any flow
of a vector field tangent to $\nf$.

From \cite{DDVG1}, in the real plane case, we deduce that for a regular value $c$ has its exponent $\rho_c = 1$
if and only if the total curvature function $t\mapsto |K|(t)$ is not continuous at $c$.

\medskip
For real polynomial with isolated singularity at infinity is $\rho_c <1$ equivalent to the continuity at
$c$ of $t \mapsto |K| (t)$ ?

More generally, we wonder if having a degenerate Gauss map, that is of rank at most $n-2$, is compatible
with condition {\bf SISI}. I really doubt it for polynomials. \\
Another way to say that is to ask, in the affine domain as well as at infinity, what conditions on the singularities
of the (generalised) critical levels (when looking at their closure in $\FP_\R^n$) of the function, a degenerate
Gauss map  is carrying ?
%
%
%
%
%
%
%
%
%
%
%
%
%
%
%
\section*{Thanks}
The author had been partially supported by the {\it European research network IHP-RAAG}
contract number HPRN-CT-2001-00271, and by {\it Deutsche Forschungs-Gemeinschaft} in the Priority
Program {\it Global Differential Geometry}.
\\
We would like to thank the University of Bath (UK) and Carl von Ossietzky Universit\"at Oldenburg (Germany)
for the working conditions provided while working on this paper.
\\
Last, but not least, a big thank you to Nicolas Dutertre for talks, questions, comments, and for a 
very careful reading that helped a lot. Cheers Mate !
%
%
%
%
%
%
%
%
%
%
%

\end{document}